\documentstyle[12pt]{article}
\def\qed{\hbox{${\vcenter{\vbox{
\hrule height 0.4pt\hbox{\vrule width 0.4pt height 6pt \kern5pt\vrule
width 0.4pt} \hrule height 0.4pt}}}$}}

\newtheorem{theorem}{Theorem}
\newtheorem{lemma}[theorem]{Lemma}
\newtheorem{definition}[theorem]{Definition} 
\newtheorem{corollary}[theorem]{Corollary} \newcommand{\proof}{\bf Proof.
\rm} \newtheorem{proposition}[theorem]{Proposition}
 
\begin{document}
%May13

\begin{titlepage}
\begin{flushright}
\hskip5cm {\it FTUV/99-40, IFIC/99-41} 
\end{flushright}

\begin{centering}

{\huge Algebra Structures on $Hom(C,L)$} 

\vspace{.5cm}
{\large G. Barnich$^1$, R. Fulp$^2$, T. Lada$^2$ and 
J. Stasheff$^{3,\dagger}$}\\

\vspace{.5cm}

$^1$ Departament de F\'{\i}sica Te\`orica, Universitat de Val\`encia,
E-46100~Burjassot (Val\`encia), Espa\~na. \\ 

\vspace{.5cm}

$^2$ Department of Mathematics, North Carolina State University, Raleigh,
NC 27695-8205.\\

\vspace{.5cm}

$^3$ Department of Mathematics, The University of North Carolina at
Chapel Hill, Phillips Hall, Chapel Hill, NC 27599-3250. 

\vspace{.5cm}

\begin{abstract}
We consider the space of linear maps from a coassociative coalgebra $C$
into a Lie algebra $L$. Unless $C$ has a cocommutative coproduct, the
usual symmetry properties of the induced bracket on $Hom(C,L)$ fail to
hold. We define the concept of twisted domain (TD) algebras in order to
recover the symmetries and also construct a modified Chevalley-Eilenberg
complex in order to define the cohomology of such algebras.
\end{abstract}

\end{centering}

\vspace{1.5cm}

{\footnotesize \hspace{-0.6cm}($^\dagger$) Research supported in part 
by NSF grant DMS-9803435.}

\end{titlepage}

\pagebreak
\section{Introduction}
The principal thrust of this note is to introduce new algebraic structures
which are defined on certain vector spaces of mappings. Assume that C is a
coassociative coalgebra and that L is an algebra which generally may be
either a Lie algebra, an associative algebra, or a Poisson algebra. We
consider the vector space Hom(C,L) of linear mappings from C into L.
When C  is cocommutative Hom(C,L) inherits whatever structure L possesses.
On the other hand when C is not cocommutative one generally obtains a new
structure on Hom(C,L) which we refer to as a TD algebra. Similar constructions
may be employed to obtain TD modules.

We are particularly interested in the case that L is a Lie algebra and in
this case we are also interested in developing a cohomology theory for
TD Lie algebras which parallels the development of Chevalley-Eilenberg
cohomology for Lie algebras. Our focus on these issues is driven by the fact
that such algebraic structures and their cohomologies arise naturally in
the study of Lagrangian field theories in physics.

More specifically, certain physical
theories, such as gauge field theories, may be formulated in terms of a
function called the
Lagrangian of the theory. From a mathematical point of view a Lagrangian is
a smooth
function defined on an appropriate jet bundle
$JE$ of some fiber bundle $E\longrightarrow M.$ If L is a Lagrangian, its
corresponding action S is defined by $S=\int_M L\  d(vol_M).$ A symmetry of the
action is a ``lift" of a generalized vector field on E to JE such that the flow
of the lifted vector field leaves S invariant. Such
symmetries are closed under the Lie brackets of generalized vector fields and
are called variational symmetries of L. A gauge symmetry of the Lagrangian is a
single member of a family of vector fields, each of which is a symmetry of S.
This family is uniquely defined by the `Noether identities' of the
Lagrangian. We interpret such a family as a mapping R from
an appropriate vector space ${\cal P}$ of "parameters" of the theory into the
space Var of variational symmetries of S. Viewed this way, if R and T are
such mappings and
$\epsilon,\eta$ are in ${\cal P},$  then $[R(\epsilon),T(\eta)]$ is a gauge
symmetry. If
$[R,T]$ is the mapping from  ${\cal P}\otimes {\cal P}$ into Var defined by
$\epsilon\otimes \eta \longrightarrow [R(\epsilon),T(\eta)],$ we see
that all possible commutators of gauge symmetries may be
realized as images of mappings from ${\cal P}\otimes {\cal P}$ into Var. In
general, if we permit mappings with domain the entire tensor coalgebra
$C=\bigoplus_n{\cal P}^{\otimes n},$ then
all possible iterated commutators of gauge symmetries are encoded as images
of such mappings. These
considerations lead
us to the study of the TD Lie structure of $Hom(C,Var).$

\section{Preliminaries}
Let $(C,\Delta )$ be a coassociative coalgebra and let $(L,\phi)$ be an algebra over a
common field $k$.  Let $Hom(C,L)$ denote the vector space of linear maps $L
\longrightarrow C$.  We  explore the properties of the convolution product
on $Hom(C,L)$ that result from  various symmetry properties of $\Delta$ and $\phi$. 
Recall that the convolution product \cite{Sweedler} on $Hom(C,L)$ is given by
$$\Phi(f,g)(c)=\phi \circ (f\otimes g)\circ \Delta(c)$$

A useful alternative to this description of $\Phi$ arises from the utilization of the
$Hom-\otimes$ interchange map \cite{MacLane} which is given by
$$\lambda :Hom(C,L)\otimes Hom(C',L')\longrightarrow Hom(C\otimes C',L\otimes
L')$$
where
$$[\lambda(f\otimes f')](c\otimes c')=f(c)\otimes f'(c').$$

Also recall that given linear maps $\alpha : C'\longrightarrow C$ and $\beta:
L\longrightarrow L'$, we have the induced maps $\alpha^*:Hom(C,L)\longrightarrow
Hom(C',L)$ and $\beta_*:Hom(C,L)\longrightarrow Hom(C,L')$.  We will make repeated use
of the equality
$$\alpha^*\circ\beta_*=\beta_*\circ\alpha^*:Hom(C,L)\longrightarrow Hom(C',L')$$
in what follows.

We thus have the convolution product described by the composition
$$Hom(C,L)\otimes Hom(C,L)\stackrel{\lambda}{\longrightarrow} Hom(C\otimes
C,L\otimes L)\stackrel{\phi_*\circ\Delta^*}{\longrightarrow}Hom(C,L),$$
i.e.
\begin{equation}
\Phi=\phi_*\circ\Delta^*\circ\lambda
\end{equation}

 We record some useful elementary properties of $\lambda$ in the
following two lemmas.
\begin{lemma}\label{naturality} (Naturality)

a)  Given a map $\zeta :A\longrightarrow C_2 $ we have the commutative diagram
$$\begin{array}{ccc}
Hom(C_1,L_1)\otimes
Hom(C_2,L_2)&\stackrel{\lambda}{\longrightarrow}&Hom(C_1\otimes C_2,L_1\otimes
L_2)\\
 & &\\
1\otimes\zeta^*\downarrow &&  \downarrow (1\otimes\zeta)^*\\
 &&\\
Hom(C_1,L_1)\otimes Hom(A,L_2) & \stackrel{\lambda}{\longrightarrow} &
Hom(C_1\otimes A,L_1\otimes L_2).
\end{array}
$$
b)  Given a map $\psi:L_2\longrightarrow B$, we have the commutative diagram
$$\begin{array}{ccc}
Hom(C_1,L_1)\otimes
Hom(C_2,L_2)&\stackrel{\lambda}{\longrightarrow}&Hom(C_1\otimes C_2,L_1\otimes
L_2)\\
 && \\
1\otimes\psi_*\downarrow &&  \downarrow (1\otimes\psi)_*\\
 && \\
Hom(C_1,L_1)\otimes Hom(C_2,B) & \stackrel{\lambda}{\longrightarrow} &
Hom(C_1\otimes C_2,L_1\otimes B).
\end{array}
$$

\end{lemma}
\proof
  For a) we have
$$[\lambda\circ(1\otimes\zeta^*)(f\otimes g)](c_1\otimes a)=[\lambda(f\otimes
g\circ\zeta)](c_1\otimes a)=f(c_1)\otimes g(\zeta(a)).$$
On the other hand,
$$[(1\otimes\zeta)^*\circ\lambda(f\otimes g)](c_1\otimes a)=\lambda(f\otimes
g)\circ(1\otimes\zeta)(c_1\otimes a)=f(c_1)\otimes g(\zeta(a)).
$$
The proof of b) follows from a similar calculation.
\qed
\medskip

It is evident that the same result holds if we replace $1\otimes\zeta^*$ by 
$\zeta^*\otimes 1$ and $1\otimes\psi_*$ by $\psi_*\otimes 1$.

\medskip
We are interested in the relationship between the symmetries of $\phi$ and
$\Delta$ and of the convolution product $\Phi$.  The next lemma will provide us
with a key link between them.  Let 
$$\Lambda : Hom(C_1,L_1)\otimes \dots \otimes Hom(C_n,L_n) \longrightarrow
Hom(C_1\otimes\dots\otimes C_n,L_1\otimes\dots\otimes L_n)$$
be given by $\Lambda=\lambda\circ(1\otimes\lambda)\circ\dots\circ(1^{n-2}\otimes\lambda)$.  Here we
are using the natural associativity of the tensor product.  When it is helpful to the
exposition, we will write $\Lambda=\Lambda^{n-1}$.
 
\begin{lemma}\label{symmetry}  (Symmetry)  Let $\sigma\in S_n$. The following diagram
commutes.
\scriptsize
$$\begin{array}{ccc}
Hom(C_1,L_1)\otimes\dots\otimes Hom(C_n,L_n) & \stackrel{\Lambda}{\longrightarrow} &
Hom(C_1\otimes\dots\otimes C_n,L_1\otimes\dots\otimes L_n) \\
 && \\
\sigma\downarrow && (\sigma^{-1})^*\circ\sigma_*\downarrow \\
 && \\
Hom(C_{\sigma(1)},L_{\sigma(1)})\otimes \dots\otimes Hom(C_{\sigma(n)},L_{\sigma(n)})
&
\stackrel{\Lambda}{\longrightarrow}&
Hom(C_{\sigma(1)}\otimes\dots\otimes C_{\sigma(n)},L_{\sigma(1)}\otimes\dots\otimes
L_{\sigma(n)})
\end{array}
$$
\normalsize
\end{lemma}
\proof  Let $f_1\otimes\dots\otimes f_n\in Hom(C_1,L_1)\otimes\dots\otimes
Hom(C_n,L_n)$.  Then by definition
$$\sigma(f_1\otimes\dots\otimes f_n)=f_{\sigma(1)}\otimes\dots\otimes
f_{\sigma(n)}$$
and further
$$[\Lambda (\sigma(f_1\otimes\dots\otimes f_n))](c_{\sigma(1)}\otimes\dots\otimes
c_{\sigma(n)})=f_{\sigma(1)}(c_{\sigma(1)})\otimes\dots\otimes
f_{\sigma(n)}(c_{\sigma(n)}).
$$
On the other hand,
$$[(\sigma^{-1})^*\circ\sigma_*\circ\Lambda(f_1\otimes\dots\otimes
f_n)](c_{\sigma(1)}\otimes\dots\otimes c_{\sigma(n)})
=[\sigma\circ\Lambda(f_1\otimes\dots\otimes
f_n)\circ\sigma^{-1}](c_{\sigma(1)}\otimes\dots\otimes c_{\sigma(n)})$$
$$=[\sigma\circ\Lambda(f_1\otimes\dots\otimes
f_n)](c_1\otimes\dots\otimes c_n)=\sigma(f_1(c_1)\otimes\dots\otimes f_n(c_n))$$
\begin{center}$=f_{\sigma(1)}(c_{\sigma(1)})\otimes\dots\otimes
f_{\sigma(n)}(c_{\sigma(n)})$. \qed
\end{center}
\medskip

As permutations will play a major role in the remainder of this note, we remark that we will
use the standard notation for cycles in the symmetric groups.  Recall that the symbol
$(i_1\,i_2\,\dots\,i_k)$ means the permutation that sends $i_1$ to $i_2$, $i_2$ to $i_3$, ...,
$i_k$ to $i_1$ ( or, if we wish, objects that are indexed by the integers $i_1$,...,$i_k$).

\section{Twisted Domain Skew Maps}

In this section we examine maps $Hom(C,L)^{\otimes n}\longrightarrow Hom(C,V)$
where $L$ and $V$ are arbitrary vector spaces. 

To begin let $(C,\Delta)$ be a coassociative coalgebra and denote the iterated coproduct
by
$$\Delta^{(n-1)}=(\Delta\otimes 1\otimes \dots\otimes 1)\circ\dots\circ
(\Delta\otimes 1)\circ\Delta :C\longrightarrow C^{\otimes n}.$$
Let $\phi:L^{\otimes n}\longrightarrow V$ be a linear map
and $\Phi:Hom(C,L)^{\otimes n}\longrightarrow Hom(C,V)$ be given by
$$\Phi = \phi_*\circ\Delta^{(n-1)*}\circ\Lambda.$$
We say that $\Phi$ is the map that is induced by the map $\phi$.
 Define also the map
$\Phi^{\sigma}:Hom(C,L)^{\otimes n}\longrightarrow Hom(C,V)$ for each $\sigma
\in S_n$ by
$$\Phi^{\sigma}=\phi_*\circ\Delta^{(n-1)*}\circ\sigma^*\circ\Lambda.$$
\medskip

The statement contained in the next lemma is a useful fact to which we will appeal
on several occasions in what follows.

\begin{lemma}\label{induced symmetry}  Let
$\Phi=\phi_*\circ\Delta^{(n-1)*}\circ\Lambda: Hom(C,L)^{\otimes n}\longrightarrow
Hom(C,V)$ be the map induced by
$\phi:L^{\otimes n}\longrightarrow V$.  Then, for $\sigma\in S_n$, the map induced
by $\phi\circ\sigma:L^{\otimes n}\longrightarrow V$ may be written as
$\Phi^{\sigma}\circ\sigma :Hom(C,L)^{\otimes n}\longrightarrow Hom(C,V)$.
\end{lemma}
\proof  The map induced by $\phi\circ\sigma$ is given by $(\phi\circ\sigma)_*
\circ\Delta^{(n-1)*}\circ\Lambda$ which is equal to
$$\phi_*\circ\sigma_*\circ\Delta^{(n-1)*}\circ\Lambda$$
$$=\phi_*\circ\Delta^{(n-1)*}\circ\sigma_*\circ\Lambda$$
$$=\phi_*\circ\Delta^{(n-1)*}\circ\sigma^*\circ(\sigma^{-1})^*\circ
\sigma_*\circ\Lambda$$
$$=\phi_*\circ\Delta^{(n-1)*}\circ\sigma^*\circ\Lambda\circ\sigma
=\Phi^{\sigma}\circ\sigma\mbox{.      \qed  }$$
\medskip

\begin{definition}  We say that the linear map $\Phi:Hom(C,L)^{\otimes n}
\longrightarrow Hom(C,V)$ is \bf twisted domain skew (TD skew) \it if
$$\Phi\circ\sigma=(-1)^{\sigma}\Phi^{\sigma^{-1}}\mbox{   for all  }\sigma\in S_n.$$
\end{definition}

Although we will use the formulation of TD skew that is given in the definition, we
may clarify the situation by considering several alternate presentations.
For $c\in C$, denote the iterated coproduct  $\Delta^{(n-1)}(c)=c_1\otimes\dots\otimes c_n$.
Then for $f_1\otimes\dots \otimes f_n\in Hom(C,L)^{\otimes n}$, we have
$$\Phi(f_1\otimes \dots\otimes f_n)(c)=\phi(f_1(c_1)\otimes\dots\otimes f_n(c_n))$$
and
$$\Phi^{\sigma}(f_1\otimes \dots\otimes f_n)(c)=\phi(f_1(c_{\sigma(1)})\otimes\dots
\otimes f_n(c_{\sigma(n)})).$$
Then the TD skew condition may be phrased as
$$\Phi(f_{\sigma(1)}\otimes\dots\otimes f_{\sigma(n)})(c)=(-1)^{\sigma}\Phi^{\sigma^{-1}}
(f_1\otimes\dots\otimes f_n)(c)$$
or
$$\phi(f_{\sigma(1)}(c_1)\otimes\dots\otimes f_{\sigma(n)}(c_n))=
(-1)^{\sigma}\phi(f_1(c_{\sigma^{-1}(1)})\otimes\dots\otimes f_n(c_{\sigma^{-1}(n)})).$$

We note here that if the coalgebra $(C,\Delta)$ is cocommutative, then the concept
of TD skew symmetry is identical to that of skew symmetry.

\medskip
\noindent
\bf Example:    \rm The example which provides the motivation for this definition is the
bracket of a Lie algebra.  Let $\phi:L\otimes L\longrightarrow L$ denote this bracket.  
The skew symmetry of $\phi$ may be expressed as $\phi\circ\tau=-\phi$ where $\tau$ is the
transposition.  The map $\Phi:Hom(C,L)\otimes Hom(C,L)\longrightarrow Hom(C,L)$ that is given
by $\Phi=\phi_*\circ\Delta^*\circ\lambda$ is not in general skew symmetric.  Indeed if we write
$\Delta(c)=x\otimes y$ and $\Phi(f\otimes g)(x\otimes y)=[f(x),g(y)]$, we
have $\Phi\circ\tau(f\otimes g)(x\otimes y)=[g(x),f(y)]=-[f(y),g(x)]=
-\Phi(f\otimes g)(y\otimes x)\neq-\Phi(f\otimes g)(x\otimes y)$.  However, in this example, we
have
$\Phi^{\tau}=\phi_*\circ\Delta^*\circ\tau^*\circ\lambda$
 and $\Phi\circ\tau=-\Phi^{\tau}$.

\medskip
 The next proposition will play a central role in the rest of this note.
\begin{proposition}\label{TDskew}  If $\phi:L^{\otimes n}\longrightarrow V$ is skew
symmetric, i.e. if
$\phi\circ\sigma=(-1)^{\sigma}\phi$ for all $\sigma\in S_n$,
 then $\Phi=\phi_*\circ\Delta^{(n-1)*}
\circ\Lambda$ is TD skew symmetric.
\end{proposition}
\proof  We calculate
$$\Phi\circ\sigma=\phi_*\circ\Delta^{(n-1)*}\circ\Lambda\circ\sigma$$
$$=\phi_*\circ\Delta^{(n-1)*}\circ(\sigma^{-1})^*\circ\sigma_*\circ\Lambda
\mbox{  (by the symmetry lemma)}$$
$$=(\phi\circ\sigma)_*\circ\Delta^{(n-1)*}\circ(\sigma^{-1})^*\circ\Lambda$$
$$=(-1)^{\sigma}\phi_*\circ\Delta^{(n-1)*}\circ(\sigma^{-1})^*\circ\Lambda
\mbox{  (because $\phi$ is skew)}$$
$$=(-1)^{\sigma}\Phi^{\sigma^{-1}}\mbox{.   \qed}$$
\medskip
Later in this note we will be confronted with various compositions of maps that
are defined on tensor products of $Hom$-sets; we would like to place such maps into the
context of maps that are induced by maps defined on the underlying vector spaces.  To be
precise, we have 
\begin {proposition}\label{comp} Suppose that we have a collection of vector spaces
$L_p$, $M_q$, $V$, and $W$ together with linear maps $$\phi:M_1\otimes\dots\otimes M_k
\longrightarrow V$$
and
$$\psi:L_1\otimes\dots\otimes L_{i-1}\otimes V\otimes L_i\otimes \dots\otimes L_n
\longrightarrow W.$$
Suppose also that $\phi$ and $\psi$ induce linear maps
$$\Phi:\bigotimes_{q=1}^k Hom(C,M_q)\longrightarrow Hom(C,V)$$
and
$$\Psi:\bigotimes_{p=1}^{i-1}Hom(C,L_p)\otimes Hom(C,V)\otimes
\bigotimes_{p=i}^n Hom(C,L_p)\longrightarrow Hom(C,W),$$
i.e., $\Phi=\phi_*\circ\Delta^{(k-1)*}\circ\Lambda^{k-1}$ and
$\Psi=\psi_*\circ\Delta^{(n)*}\circ\Lambda^n$.  Then the composition
$$\Psi\circ(1^{\otimes(i-1)}\otimes\Phi\otimes 1^{\otimes(n-i+1)}):$$
$$\bigotimes_{p=1}^{i-1}Hom(C,L_p)\otimes\bigotimes_{q=1}^k Hom(C,M_q)\otimes
\bigotimes_{p=i}^n Hom(C,L_p)\longrightarrow Hom(C,W)$$
is induced by $\psi\circ(1^{\otimes(i-1)}\otimes\phi\otimes 1^{\otimes(n-i+1)})$, i.e.
$$\Psi\circ(1^{\otimes(i-1)}\otimes\Phi\otimes 1^{\otimes(n-i+1)})=
[\psi\circ(1^{\otimes(i-1)}\otimes\phi\otimes 1^{\otimes(n-i+1)})]_*
\circ\Delta^{(n+k-1)*}\circ\Lambda^{n+k-1}.$$
\end{proposition}
\proof
The coassociativity of $\Delta$ gives us the equality
$$\Delta^{(n+k-1)}=(1^{\otimes(i-1)}\otimes\Delta^{(k-1)}\otimes 1^{(n-i+1)})
\circ\Delta^{(n)}$$
and the associativity of the tensor product  gives us
$$\Lambda^{n+k-1}=\Lambda^n\circ(1^{\otimes(i-1)}\otimes\Lambda^{k-1}
\otimes 1^{\otimes(n-i+1)}).$$
As a result, we have
$$[\psi\circ(1^{\otimes(i-1)}\otimes\phi\otimes 1^{\otimes(n-i+1)})]_*
\circ\Delta^{(n+k-1)*}\circ\Lambda^{n+k-1}$$
\medskip
$$=\psi_*\circ(1^{\otimes(i-1)}\otimes\phi\otimes 1^{\otimes(n-i+1)})_*
\circ((1^{\otimes(i-1)}\otimes\Delta^{(k-1)}\otimes 1^{\otimes(n-i+1)})
\circ\Delta^{(n)})^*$$
$$\circ\Lambda^n\circ(1^{\otimes(i-1)}\otimes\Lambda^{k-1}
\otimes1^{\otimes(n-i+1)})$$
\medskip
$$=\psi_*\circ\Delta^{(n)*}\circ(1^{\otimes(i-1)}\otimes\phi
\otimes 1^{\otimes(n-i+1)})_*
\circ((1^{\otimes(i-1)}\otimes\Delta^{(k-1)}\otimes 1^{\otimes(n-i+1)})^*$$
$$\circ\Lambda^n\circ(1^{\otimes(i-1)}\otimes\Lambda^{k-1}
\otimes1^{\otimes(n-i+1)})$$
\medskip
$$=\psi_*\circ\Delta^{(n)*}\circ\Lambda^n\circ(1^{\otimes(i-1)}\otimes\phi_*
\otimes 1^{\otimes(n-i+1)})
\circ((1^{\otimes(i-1)}\otimes\Delta^{(k-1)*}\otimes 1^{\otimes(n-i+1)})$$
$$\circ(1^{\otimes(i-1)}\otimes\Lambda^{k-1}
\otimes1^{\otimes(n-i+1)})$$
(by several applications of Lemma \ref{naturality})
\medskip
$$=(\psi_*\circ\Delta^{(n)*}\circ\Lambda^n)\circ(1^{\otimes(i-1)}\otimes
\phi_*\circ\Delta^{(k-1)*}\circ\Lambda^{k-1}\otimes 1^{\otimes(n-i+1)})$$
\medskip
\begin{center}
$=\Psi\circ(1^{\otimes(i-1)}\otimes\Phi\circ 1^{\otimes(n-i+1)})$.  \qed
\end{center}
\medskip
If we have $\sigma\in S_{n+k}$, we will denote the map $$[\psi\circ
(1^{\otimes(i-1)}\otimes\phi\otimes 1^{(n-i+1)})]_*\circ\Delta^{(n+k-1)*}
\circ\sigma^*\circ\Lambda^{n+k-1}$$ by $(\Psi\circ(1^{\otimes(i-1)}\otimes
\Phi\circ 1^{\otimes(n-i+1)}))^{\sigma}$, or by
$$(\Psi\circ(1^{\otimes(i-1)}\otimes
\Phi\circ 1^{\otimes(n-i+1)}))^{\sigma}(g_1,\dots,g_{n+k})$$
$$= (\Psi(g_1,\dots,g_{i-1},\Phi(g_i,\dots,g_{i+k-1}),g_{i+k},
\dots,g_{n+k}))^{\sigma}.$$

\medskip
   
\section{TD Algebras}
When the product on an algebra $L$ satisfies some (skew) symmetry relations, the induced
product on $Hom(C,L)$ will satisfy the corresponding twisted domain symmetries.  Indeed
if $(L,\phi)$ is an algebra with relations given by permutations, i.e.
$F(\phi_{\alpha}^i,\phi_{\alpha}^i\circ\sigma_{\beta})=0$, then the corresponding TD algebra
will have relations
$F(\Phi_{\alpha}^i,(\Phi_{\alpha}^i)^{\sigma_{\beta}}\circ\sigma_{\beta})=0$ where the
$\phi_{\alpha}^i$'s  are various iterates of $\phi$ and $1\otimes\dots\phi\dots\otimes 1$.

We also remark that relations in an algebra that do not involve permutations carry over intact
to the convolution product.  The fundamental example, of course, is associativity; i.e. if the
algebra $(L,\phi)$ is associative and the coalgebra $(C,\Delta)$ is coassociative, then  the
algebra $(Hom(C,L),\Phi)$ is associative.

In this section we examine these TD symmetries for several
types of algebras.  Again, we assume that
$(C,\Delta )$ is a fixed coassociative coalgebra.  We are primarily interested in the example
where $(L,
\phi)$ is a Lie algebra.  Recall that the basic relations
 in a Lie algebra may be phrased as follows:
\begin{equation}
\mbox{skew symmetry:  }\phi\circ\tau=-\phi
\end{equation}
and
\begin{equation}
\mbox{Jacobi identity:   } \phi\circ(1\otimes\phi)+\phi\circ(1\otimes\phi)\circ\xi+
\phi\circ(1\otimes\phi)\circ\xi^2=0
\end{equation}
where $\tau:L\otimes L\longrightarrow L\otimes L$ is the transposition and
$\xi:L\otimes L\otimes L\longrightarrow L\otimes L\otimes L$ is the cyclic
permutation (123).  
\begin{proposition}  Let ($L,\phi)$ be a Lie algebra and let
$\Phi=\phi_*\circ\Delta^*\circ\lambda$. Then
$(Hom(C,L),\Phi)$ has the following symmetry properties:

\begin{equation}
\Phi\circ\tau=-\Phi^{\tau}
\label{skew}
\end{equation}
and
\begin{equation}
\Phi\circ(1\otimes\Phi)+\Phi\circ(1\otimes\Phi)^{\xi}\circ\xi
+\Phi\circ(1\otimes\Phi)^{\xi^2}\circ\xi^2=0
\label{Jacobi}
\end{equation}
\end{proposition}
\proof  
Equation \ref{skew} follows directly from Proposition \ref{TDskew}. Next we 
rewrite the left hand side of Equation \ref{Jacobi} as  
\begin{eqnarray}
\phi_*\circ(1\otimes\phi)_*\circ((1\otimes\Delta)\circ\Delta)^*\circ\Lambda
\nonumber\\
+\phi_*\circ(1\otimes\phi)_*\circ(\xi\circ(1\otimes\Delta)\circ\Delta)^*
\circ\Lambda\circ\xi\nonumber\\
+\phi_*\circ(1\otimes\phi)_*\circ(\xi^2\circ(1\otimes\Delta)\circ\Delta)^*
\circ\Lambda\circ\xi^2
\nonumber
\end{eqnarray}
\medskip
$$=\phi_*\circ(1\otimes\phi)_*\circ((1\otimes\Delta)\circ\Delta)^*\circ\Lambda$$
$$+\phi_*\circ(1\otimes\phi)_*\circ(\xi\circ(1\otimes\Delta)\circ\Delta)^*
\circ(\xi^2)^*\circ\xi_*\circ\Lambda$$
$$+\phi_*\circ(1\otimes\phi)_*\circ(\xi^2\circ(1\otimes\Delta)\circ\Delta)^*
\circ\xi^*\circ(\xi^2)_*\circ\Lambda$$
by Lemma \ref{symmetry}, and 
\medskip
$$=\phi_*\circ(1\otimes\phi)_*\circ((1\otimes\Delta)\circ\Delta)^*\circ\Lambda$$
$$+\phi_*\circ(1\otimes\phi)_*\circ\xi_*\circ(\xi^2\circ\xi\circ(1\otimes\Delta)
\circ\Delta)^*
\circ\Lambda$$
$$+\phi_*\circ(1\otimes\phi)_*\circ\xi^2_*\circ(\xi\circ\xi^2\circ(1\otimes
\Delta)\circ\Delta)^*\circ\Lambda$$

\medskip
$$=\phi_*\circ(1\otimes\phi)_*\circ(1+\xi_*+\xi^2_*)\circ((1\otimes\Delta)
\circ\Delta)^*\circ\Lambda$$
$$=[\phi\circ(1\otimes\phi)\circ(1+\xi+\xi^2)]_*\circ((1\otimes
\Delta)\circ\Delta)^*\circ\Lambda=0$$
because $\phi$ satisfies the Jacobi identity.  We also used the fact that 
$\xi^{-1}=\xi^2$
and that the iterated convolution product
$$\Phi\circ(1\otimes\Phi)=
\phi_*\circ(1\otimes\phi)_*\circ\Delta^*\circ(1\otimes\Delta)^*\circ
\Lambda$$
by Proposition \ref{comp}.  \qed

\medskip
When the coalgebra $(C,\Delta)$ is cocommutative, the usual Lie algebra structure on
$Hom(C,L)$ is evident.
\begin{corollary}  If the coassociative coalgebra $(C,\Delta)$ is cocommutative, i.e.
$\Delta=\tau\circ\Delta$, then $\Phi$ is a Lie algebra structure for $Hom(C,L)$.
\end{corollary}
\proof  The cocommutativity of $\Delta$ implies that 
$\xi\circ(1\otimes\Delta)\circ\Delta =(1\otimes\Delta)\circ\Delta$.  We may write
$\xi=(\tau\otimes 1)\circ(1\otimes\tau)$ and then $\xi\circ(1\otimes\Delta)\circ\Delta =
(\tau\otimes 1)\circ(1\otimes\tau)\circ(1\otimes\Delta)\circ\Delta = (\tau\otimes
1)\circ(1\otimes\Delta)\circ\Delta = (\tau\otimes 1)\circ(\Delta\otimes 1)\circ\Delta=
(\Delta\otimes 1)\circ\Delta=(1\otimes\Delta)\circ\Delta$.  Replacement of $\tau\circ\Delta$ by
 $\Delta$ in equation \ref{skew} yields the skew
symmetry relation;  replacement of $\xi\circ(1\otimes\Delta)\circ\Delta$ and
$\xi^2\circ(1\otimes\Delta)\circ\Delta$ by $(1\otimes\Delta)\circ\Delta$ yields the usual
Jacobi identity.  \qed

\medskip
Another interesting example is the following:
\begin{corollary} \label{jordan}  If the coassociative coalgebra $(C,\Delta)$ is skew
cocommutative, i.e.
$\tau\circ\Delta=-\Delta$, then the bracket $\Phi$ is symmetric and satisfies the Jacobi
identity.
\end{corollary}
\proof  From the proof of Proposition \ref{TDskew}, we have
$$\Phi \circ\tau=(\phi\circ\tau)_*\circ(\tau\circ\Delta)^*\circ\lambda$$
$$=-\phi_*\circ(\tau\circ\Delta)^*\circ\lambda=\phi_*\circ\Delta^*\circ\lambda$$

For the Jacobi identity we have almost as in the previous proof
$$\xi\circ(1\otimes\Delta)\circ\Delta =(\tau\otimes
1)\circ(1\otimes\tau)\circ(1\otimes\Delta)\circ\Delta$$
$$ = -(\tau\otimes
1)\circ(1\otimes\Delta)\circ\Delta =- (\tau\otimes 1)\circ(\Delta\otimes 1)\circ\Delta$$
$$ = 
(\Delta\otimes 1)\circ\Delta=(1\otimes\Delta)\circ\Delta.$$
A similar calculation holds for $\xi^2\circ(1\otimes\Delta)\circ\Delta$.  \qed

\medskip
\noindent
\bf Remark:    \rm Let us write $\Phi(f\otimes g)=fg$ and then the conditions in Corollary
\ref{jordan} imply the equalities $fg=gf$ and $(f^2g)f+f^2(gf)=0$.  These are the defining
equations of a
\bf Jordan algebra\rm.  To see that the second equation holds, we write the Jacobi identity as
$(fg)h+(hf)g+(gh)f=0$;  if we first let $g=h=f$ we obtain the equality$(f^2)f+(f^2)f+(f^2)f=0$,
i.e.$(f^2)f=0$ when characteristic $k\neq 3$.  Next replace $h$ by $f$ and $f$ by $f^2$ in the
 Jacobi
identity and obtain the equation $(f^2g)f+(ff^2)g+(gf)f^2=0$.  The middle term drops out and the
remainder (after several applications of commutativity) yields the result.

\medskip
Of course, the coalgebras that we have in mind for the three previous results are the
tensor coalgebra $T^CV$ generated by the vector space $V$ and the subcoalgebras of
symmetric tensors and skew symmetric tensors.  

\medskip

We conclude this section with one more example.
Recall that the vector space $L$ is a \bf Poisson algebra \rm if it possesses  a Lie bracket
$\phi$ as well as an associative, commutative multiplication $\mu$;  also required is the
property that the bracket $\phi$ acts as a derivation with respect to the operation $\mu$.

This derivation property is usually written as
$$[a,bc]=[a,b]c+b[a,c]$$
where $\phi(a,b)=[a,b]$ and $\mu(a,b)=ab$.  We will however use the commutativity of $\mu$ to
rewrite this relation as
\begin{equation}
[a,bc]=c[a,b]+b[a,c]
\end{equation}
The reason for this is that we can better describe this relation in functional notation
together with permutations as
\begin{equation}
\phi\circ(1\otimes\mu)=\mu\circ(1\otimes\phi)\circ\xi+\mu\circ(1\otimes\phi)
\circ(\tau\otimes
1)
\end{equation}
where as before $\xi$ is the permutation on $L\otimes L\otimes L$ given by (123) and
$\tau\otimes 1$ is the permutation given by (12).

\begin{proposition}  Let $(L,\phi,\mu)$ be a Poisson algebra and let $(C,\Delta)$ be a
coassociative coalgebra.  Let $\Phi=\phi_*\circ\Delta^*\circ\lambda$ and $M
=\mu_*\circ\Delta^*\circ\lambda$.  Then
$(Hom(C,L),\Phi,M)$ has the structure of a TD Poisson algebra;  i.e. $\Phi$
satisfies the Lie relations up to permutation as in Proposition 3, $M$ is commutative
up to permutation, and the derivation property holds up to permutation.
\end{proposition}
\proof  The TD Lie structure for $\Phi$ was shown in Proposition 3. For the
TD commutativity we have
$$M\circ\tau=\mu_*\circ\Delta^*\circ\lambda\circ\tau=\mu_*\circ\Delta^*\circ
\tau_*\circ\tau^*\circ\lambda$$
$$=(\mu\circ
\tau)_*\circ(\tau\circ\Delta)^*\circ\lambda=\mu_*\circ(\tau\circ
\Delta)^*\circ\lambda=M^{\tau}.$$

For the derivation property we have
$$\Phi\circ(1\otimes M)=\phi_*\circ(1\otimes \mu)_*\circ\Delta^*
\circ(1\otimes\Delta)^*\circ\Lambda$$
$$=(\mu\circ(1\otimes\phi)\circ\xi)_*\circ\Delta^*\circ(1\otimes
\Delta)^*\circ\Lambda
+(\mu\circ(1\otimes\phi)\circ(\tau\otimes 1))_*\circ\Delta^*\circ(1\otimes
\Delta)^*\circ\Lambda$$
(because $(L,\phi,\mu)$ is a Poisson algebra)
$$=(\mu\circ(1\otimes\phi))_*\circ\Delta^*\circ(1\otimes\Delta)^*\circ\xi_*
\circ\Lambda$$
$$+(\mu\circ(1\otimes\phi))_*\circ\Delta^*\circ(1\otimes\Delta)^*
\circ(\tau\otimes1)_*\circ\Lambda$$
\medskip
$$=(\mu\circ(1\otimes\phi))_*\circ\Delta^*\circ(1\otimes\Delta)^*\circ\xi^*
\circ\Lambda\circ\xi$$
$$+(\mu\circ(1\otimes\phi))_*\circ\Delta^*\circ(1\otimes\Delta)^*\circ
(\tau\otimes
1)^*\circ\Lambda\circ (\tau\otimes 1)$$
which is the derivation property except for the presence of the permutations 
$\xi^*$ and
$(\tau\otimes 1)^*$ in the final two lines.  \qed

These examples will be studied further in future work.

\section{TD Module Structures and Cohomology}
\subsection{TD modules}
Let $(L,\phi)$ be a Lie algebra and let $B$ be a module over $L$ via the linear
map $\psi:L\otimes B\longrightarrow B$ which satisfies the equality
$$\psi(\phi(x_1,x_2),b)=\psi(x_1,\psi(x_2,b))-\psi(x_2,\psi(x_1,b))$$
where $x_1,x_2\in L\mbox{ and }b\in B$.  Let us write this relation as
$$\psi\circ(\phi\otimes 1)=\psi\circ(1\otimes \psi)-\psi\circ(1\otimes\psi)\circ\tau$$
where $\tau=(12)\in S_3$.  We may use this formulation to make the following definition.

\begin{definition}  Let $(Hom(C,L) ,\Phi)$ be a TD Lie algebra.  Then the vector
space $Hom(C,B)$ is a
 module over $Hom(C,L)$ if there is a linear map
$$\Psi:Hom(C,L)\otimes Hom(C,B)\longrightarrow Hom(C,B) $$
induced by $\psi:L\otimes B\longrightarrow B$ that satisfies the equation
$$\Psi\circ(\Phi\otimes 1)=\Psi\circ(1\otimes \Psi)-(\Psi\circ(1\otimes\Psi))^{\tau}
\circ\tau.$$
\end{definition}
\medskip
The next proposition follows immediately from the above definition.

\begin{proposition}  Let $\psi:L\otimes B\longrightarrow B$ give $B$ the structure of
a Lie module over $L$.  Then $\Psi=\psi_*\circ\Delta^*\circ\lambda$ gives 
$Hom(C,B)$ the structure of a $Hom(C,L)$ module.
\end{proposition}
\medskip
\noindent
\bf Example:  \rm The fundamental example (indeed the example that motivates the
definition) of a TD Lie module is the following:  let $(Hom(C,L),\Phi)$ be a
TD Lie algebra.  Then $Hom(C,L)$ is a module over itself with structure map
given by $\Phi$.  The proof of this is exactly the same as the Lie algebra case
except for the twisting in the Jacobi identity.

\medskip
\subsection{Cohomology}
We next consider the cohomology of a
 TD Lie algebra with coefficients in a module over it.  We begin by reviewing
 the Chevalley - Eilenberg complex for Lie algebras.  Let $(L,\phi)$ be a Lie algebra and
$B$ a Lie module over
$L$ with structure map $\psi$.  Let $Alt^n(L,B)$ denote the vector space of skew-
symmetric linear maps $L^{\otimes n}\longrightarrow B$.  Define a linear map
$$d:Alt^n(L,B)\longrightarrow Alt^{n+1}(L,B)$$
for $f_n\in Alt^n(L,B)$ by
$$df_n(x_1,\dots,x_{n+1})=\sum^{n+1}_{i=1}(-1)^{i+1}\psi(x_i,f_n(x_1,\dots,
\hat{x_i},\dots,x_{n+1}))$$
$$+\sum_{j<k}(-1)^{j+k}f_n(\phi(x_j,x_k),x_1,\dots,\hat{x_j},\dots,\hat{x_k},
\dots,x_{n+1}).$$
One can check that $d^2=0$ and thus we have a cochain complex.  A useful fact about
this differential is that each of the two summands is a skew symmetric map.  This may be
verified by using the next two lemmas.

\begin{lemma}\label{skew1}  Suppose that $f_n:L^{\otimes n}\longrightarrow V$ is a skew
symmetric map.  Then the extension of $f_n$ to the map $f:L^{\otimes n+1}\longrightarrow V$ 
that is given by $f(x_1,\dots,x_{n+1})=\sum_{i=1}^n(-1)^{i+1}x_i\otimes
f_n(x_1,\dots,\hat{x_i},\dots,x_{n+1})$ is skew symmetric.
\end{lemma}
\proof We verify the claim for the transposition $(p\,p+1)$.
$$f(x_1,\dots,x_{p+1},x_p,\dots,x_{n+1})=\sum_{i\neq p,p+1}(-1)^{i+1}x_i\otimes
f_n(x_1,\dots,\hat{x_i},\dots,x_{p+1},x_p,\dots,x_{n+1})$$
$$+(-1)^{p+1}x_{p+1}\otimes f_n(x_1,\dots,\hat{x_{p+1}},x_p,\dots,x_{n+1})\;\;\mbox{     
  ($x_{p+1}$ is the $p$-th coordinate)}$$
$$+(-1)^{p+2}x_p\otimes f_n(x_1,\dots,x_{p+1},\hat{x_p},\dots,x_{n+1})\;\;\mbox{     ($x_p$ is
the $p+1$-st coordinate)}$$
\medskip
$$=-\sum_{i\neq p,p+1}(-1)^{i+1}x_i\otimes
f_n(x_1,\dots,\hat{x_i},\dots,x_p,x_{p+1},\dots,x_{n+1})$$
$$-(-1)^{p+2}x_{p+1}\otimes f_n(x_1,\dots,\hat{x}_{p+1},\dots,x_{n+1})$$
$$-(-1)^{p+1}x_p\otimes f_n(x_1,\dots,\hat{x_p},\dots,x_{n+1})$$
\medskip
$$=-f(x_1,\dots,x_p,x_{p+1},\dots,x_{n+1}).\mbox{   \qed}$$ 

\medskip
\begin{lemma}\label{skew2}  Suppose that $f_k:L^{\otimes k}\longrightarrow L$ and
$f_{n-k+1}:L^{\otimes n-k+1}\longrightarrow L$ are skew symmetric maps. Then the map
$f:L^{\otimes n}\longrightarrow L$ given by 
$$f(x_1,\dots,x_n)=\sum_{\sigma}(-1)^{\sigma}f_{n-k+1}(f_k(x_{\sigma(1)},
\dots,x_{\sigma(k)}),\dots,x_{\sigma(n)})$$
where $\sigma$ runs through all $(k,\,n-k)$ unshuffles, is skew symmetric.
\end{lemma}
\proof  Again, we verify the claim for the transposition $(p\,p+1)$ and show that
$$f(x_1,\dots,x_{p+1},x_p,\dots,x_n)=-f(x_1,\dots,x_p,x_{p+1},\dots,x_n).$$
In the expansion of $f$, there are three situations to consider.  In the first, the indices
$p,\;p+1\in\{\sigma(1),\dots\sigma(k)\}$ and in the second,
$p,\,p+1\in\{\sigma(k+1),\dots,\sigma(n)\}$. The skew symmetry of $f_k$ takes care of the
first case while the skew symmetry of $f_{n-k+1}$ takes care of the second.  The remaining
case occurs when $p$ and $p+1$ are in different sets that are given by the unshuffle
decomposition.  However, such a term pairs off with the negative of the corresponding term in
the expansion of $f(x_1,\dots,x_p,x_{p+1},\dots,x_n)$ because the permutations that lead to
these terms differ only by the product with the transposition $(p\,p+1)$.   \qed 
\medskip

The first summand in the differential is the composition of a linear map, $\psi$, with
the extension of the map $f_n$ using $(1,n)$ - unshuffles; the second is the composition
of the skew symmetric map $f_n$ with the extension of the bracket, $\phi$, using
$(2,n-1)$ - unshuffles.  Consequently, we have the alternate description of the differential in
the Chevalley-Eilenberg complex given by
$$df=\sum_{\sigma}(-1)^{\sigma}\psi\circ(1\otimes f)\circ\sigma -
\sum_{\sigma'}(-1)^{\sigma'}f\circ(\phi\otimes 1)\circ\sigma'$$
where $\sigma$ runs through all $(1,n)$ - unshuffles and $\sigma'$ runs through all
$(2,n-1)$ - unshuffles.
\medskip

We now construct a cochain complex that may be used to calculate the cohomology of the
TD Lie algebra $(Hom(C,L),\Phi)$ with coefficients in the module $Hom(C,B)$ with
structure map $\Psi$.  Let $TDalt^n(Hom(C,L),Hom(C,B))$ be the vector space of TD
skew symmetric maps $(Hom(C,L))^{\otimes n}\longrightarrow Hom(C,B)$. For $n=0$,
$TDalt^0(Hom(C,L),Hom(C,B))$ consists of constant maps from $Hom(C,L)$ to the constant maps
from $C$ to $B$.  Moreover, these maps must be induced by constant maps from $L$ to $B$.
Consequently, we define
$TDalt^0$ to be the vector space $B$.  We also define $TDalt^1(Hom(C,L),Hom(C,B))$ to be the
vector space of linear maps that are induced by linear maps $L\longrightarrow B$

 Define a linear map
$$\delta:TDalt^n(Hom(C,L),Hom(C,B))\longrightarrow TDalt^{n+1}(Hom(C,L),Hom(C,B))$$
as follows: for $F_n\in TDalt^n$ that is induced from $f_n\in Alt^n$,
$$\delta F_n(g_1,\dots\, g_{n+1})=\sum^{n+1}_{i=1}(-1)^{i+1}\Psi(g_i,
F_n(g_1,\dots,\hat{g_i},\dots,g_{n+1}))^{\sigma}$$
$$-\sum_{j<k}(-1)^{j+k}F_n(\Phi(g_j,g_k),g_1,\dots,\hat{g_j},\dots,\hat{g_k},
\dots,  g_{n+1})^{\sigma^{\prime}}$$
where $\sigma$ is the $(1,n)$ unshuffle that in effect interchanges the first and the $i$-th
coordinates and
$\sigma'$ is the $(2,n-1)$ unshuffle that replaces the first and second coordinates by the
$j$-th and $k$-th coordinates.

When $n=0$, we have $(\delta f_0)(g)(x)=\Psi(g,f_0)(x)=\psi(g(x),f_0)$ for
$g\in Hom(C,L),x\in C$.  Here we used the identification $F_0=f_0\in B$.  

Let us write $\delta$ in the form
$$\delta F=\sum_{\sigma}(-1)^{\sigma}(\Psi\circ(1\otimes
F))^{\sigma}\circ\sigma-\sum_{\sigma'}(-1)^{\sigma'}(F\circ(\Phi\otimes 1))^{\sigma'}\circ
\sigma'.$$

\begin{proposition} $\delta$ is well defined and is induced by $d$.
\end{proposition}
\proof We first show that $\delta$ is induced by $d$.  Begin by considering the case in which
both $\sigma$ and $\sigma'$ equal the identity element in $S_n$.  Then the map induced by
$\psi\circ(1\otimes f)$, for $\psi\in Hom(L\otimes B,B)$ and $f\in 
Hom(L^{\otimes n},B)$, is the map $\Psi\circ(1\otimes F)$ by Proposition \ref{comp}.  Here, we assume
that $F$ is induced by $f$ and $\Psi$ is induced by $\psi$.  Similarly,
 the map induced by $f\circ(\phi\otimes 1)$ is the map $F\circ(\Phi\otimes 1)$.

Now for general $\sigma$, $\sigma'$, the map $\psi\circ(1\otimes f)\circ\sigma$ induces the map
$\Psi\circ(1\otimes F)^{\sigma}\circ\sigma$ and the map $f\circ(\phi\otimes
1)\circ\sigma'$ induces the map $F\circ(\Phi\otimes 1)^{\sigma'}\circ\sigma'$ by
Lemma \ref{induced symmetry}.  Because the maps $\sum_{\sigma}(-1)^{\sigma}\psi
\circ(1\otimes f)\circ\sigma$ and $\sum_{\sigma'}(-1)^{\sigma'}f\otimes
(\phi\otimes 1)\sigma'$ are skew symmetric, the induced maps $\sum_{\sigma}(-1)^
{\sigma}\Psi\circ(1\otimes F)^{\sigma}\circ\sigma$ and $\sum_{\sigma'}(-1)^
{\sigma'}F\circ(\Phi\otimes 1)^{\sigma'}\circ\sigma'$ are TD skew symmetric by 
Proposition \ref{TDskew} and it follows that $\delta$ is well defined.
\qed
\medskip
\begin{corollary} $\delta^2=0$
\end{corollary}
\proof  This now follows directly from the fact that $d^2=0$.  \qed

Of course, we now define the cohomology of $Hom(C,L)$ to be the cohomology of this complex.
Just as in the Lie algebra setting, \newline
$H^0(Hom(C,L),Hom(C,B))$ will equal the subspace
of invariant elements $\beta\in Hom(C,B)=B$ in the sense that $\Psi(\alpha,\beta) =0$ for all
$\alpha\in Hom(C,L)$.
\bigskip

\section{Twisted Domain Lie - Rinehart}
In the classical version of Lie Rinehart cohomology \cite{R}, the situation involves a Lie
algebra $(L,\phi )$ and an associative algebra $(B,\nu)$ which are interrelated by
the following data.  First of all, $L$ is a left $B$-module via a structure map
$\mu:B\otimes L\longrightarrow L$; we usually write $\mu(a,x)=ax$.  Next we have that
$B$ is an $L$ module via a structure map $\psi:L\otimes B\longrightarrow B$;
moreover, we assume that for fixed $x\in L$, $\psi$ is a derivation on $B$; i.e.
$\psi(x,ab)=a\psi(x,b)+\psi(x,a)b$ where we write $\nu(a,b)=ab$.
  Recall that this means that the adjoint
$\psi^*$  of $\psi$ is a Lie homomorphism from $L$ to the Lie algebra of derivations on $B$.

These structures are further related by the following two conditions:

\begin{itemize}
\item [LRa:]     $\psi\circ(\mu\otimes 1)=\nu\circ(1\otimes \psi)$
\item [LRb:]     $\phi\circ(1\otimes\mu)=\mu\circ[(1\otimes\phi)\circ\tau
+(\psi\otimes 1)]$
\end{itemize}
where $\tau=(1\,2)\in S_3$.
A pair $(L,B)$ that satisfies all of the above conditions is called a \bf Lie - 
Rinehart pair\rm.

The property labeled LRa may be rephrased as $\psi(ax,b)=a\psi(x,b)$ for $x\in L$
and $a,b\in B$ and says that the map $\psi$ is $B$-linear.

The property labeled LRb may be rephrased as $[x,ay]=\psi(x,a)y+a[x,y]$ for
$a\in B$ and $x,y\in L$ and describes the relationship between the bracket on $L$
with the two module structures.

Such structures are of interest when one wishes to study the maps in $Alt^n(L,B)$
that are $B$-linear.  The conditions required for a Lie - Rinehart pair will guarantee 
that the subspace of $Alt^n(L,B)$ of $B$-linear maps together with the restriction of
the Chevalley-Eilenberg differential is in fact a subcomplex of the Chevalley-Eilenberg
complex.

We next examine the structures on the pair $(Hom(C,L),Hom(C,B))$ that are induced  by a Lie-Rinehart
structure on $(L,B)$.

As in the previous section we have the twisted Lie algebra $(Hom(C,L),\Phi)$ and the module
$Hom(C,B)$ with structure map $\Psi$.  The associative algebra structure on $(B,\nu)$ induces
an associative (no twist) algebra structure on $(Hom(C,B),\bar{\nu})$.  The left
$B$-module structure on $L$ induces a left $Hom(C,B)$-module structure on
$Hom(C,L)$ (again, no twist needed)
$$\bar{\mu}:Hom(C,B)\otimes Hom(C,L)\longrightarrow Hom(C,L)$$
such that $\bar{\mu}\circ(1\otimes\bar{\mu})=\bar{\mu}\circ(\bar{\nu}\otimes
1)$.  Moreover, $\Psi$ is a twisted derivation on
$Hom(C,B)$; i.e.
$$\Psi\circ(1\otimes \bar{\nu})=(\bar{\nu}\circ(\Psi\otimes
1))^{\tau}\circ\tau+\bar{\nu}\circ(\Psi\otimes 1)$$
where $\tau=(1\,2)\in S_3$, or, for fixed $g\in Hom(C,L)$,
$$\Psi(g,\bar{\nu}(\alpha,\beta))=\bar{\nu}(\alpha,\Psi(g,\beta))^{\tau}
+\bar{\nu}(\Psi(g,\alpha),\beta).$$

The appropriate requirements for twisted Lie-Rinehart now take the form

\bigskip
\noindent
TDLRa:   $\Psi\circ(\bar{\mu}\otimes 1)=\bar{\nu}\circ(1\otimes\Psi)$
\newline
\noindent
TDLRb:   $\Phi\circ(1\otimes\bar{\mu})=[\bar{\mu}\circ(1\otimes\Phi)]^{\tau}\circ
\tau+\bar{\mu}\circ(\Psi\otimes 1)$

\bigskip
\noindent
where $\tau=(1\,2)
\in S_3$.
We say that $(Hom(C,L),Hom(C,B))$ is a \bf TD Lie Rinehart \rm pair.  

Observe that by using skew symmetry, condition LRb may be written in the form 
$$\phi\circ(\mu\otimes 1)=\mu\circ(1\otimes\phi)-\mu\circ(\psi\otimes 1)\circ\sigma$$
where $\sigma=(1\,2\,3)\in S_3$.  The corresponding TDLRb thus assumes the form
$$\Phi\circ(\bar{\mu}\otimes 1)=\bar{\mu}\circ(1\otimes\Phi)-\bar{\mu}\circ(\Psi\otimes
1)^{\sigma}\circ\sigma$$
which in turn may be written as
$$\Phi(\bar{\mu}(\beta,g_1),g_2)=\bar{\mu}(\beta,\Phi(g_1,g_2))-
\bar{\mu}(\Psi(g_2,\beta),g_1)^{\sigma}$$
where $g_1,g_2\in Hom(C,L),\beta\in Hom(C,B)$ and $\sigma=(1\,2\,3)\in S_3$.

We next consider elements $\alpha$ of $TDalt^n(Hom(C,L),Hom(C,B))$ that are
  ${\bf Hom(C,B)}${\bf-linear}.  By this we mean that
$$\alpha\circ(1^{\otimes(i-1)}\otimes\bar{\mu}\otimes 1^{\otimes(n-i+1)})
(g_1,\dots,g_{i-1},\beta,g_i,\dots,g_n)=(\bar{\nu}\circ(1\otimes\alpha))
^{\sigma}(\beta,g_1,\dots,g_n)$$
or
$$\alpha(g_1,\dots,g_{i-1},\bar{\mu}(\beta,g_i),g_{i+1},\dots,g_n)=
\bar{\nu}(\beta,\alpha(g_1,\dots,g_n))^{\sigma}$$
for each $i$, where $\sigma$ is the cyclic permutation $(0\,1\dots(i-1))$ 
in $S_{n+1}$ regarded as acting on the ordered set $\{0,1,\dots,n\}$, $\beta\in
Hom(C,B)$, and
$(g_1,\dots,g_n)
\in Hom(C,L)^{\otimes n}$.

The point here is that just as in the Lie algebra case where the Lie Rinehart complex is a subcomplex 
of the Chevalley-Eilenberg complex, the TD Lie-Rinehart complex is a subcomplex of the TD
Chevalley-Eilenberg complex.

To see this, we need only verify that the Chevalley-Eilenberg differential $\delta$ preserves
$Hom(C,B)$ linearity.  We summarize this in 
\begin{theorem}If $(Hom(C,L),Hom(C,B))$ is a TD Lie-Rinehart pair, then the TD Lie-Rinehart
complex is a subcomplex of the TD Chevalley-Eilenberg complex.
\end{theorem}
\proof  We claim that 
$$\delta F_n(g_1,\dots,\bar{\mu}(\beta,g_i),\dots,g_{n+1})=
\bar{\nu}(\beta,\delta F_n(g_1,\dots,g_i,\dots,g_{n+1}))^{\xi}$$
for each $\xi=(0\,1\dots i-1)\in S_{n+2}$.  Here, we regard the permutations as acting on the 
set of integers $\{0,1,\dots,n+1\}$.
We first verify this for the case $i=1$ and then use skew symmetry to complete the proof.

To begin, we must show 
that
$$\delta F_n(\bar{\mu}(\beta,g_1),g_2,\dots,g_{n+1})=\bar{\nu}(\beta,\delta
F_n(g_1,\dots,g_i,\dots,g_{n+1})).$$
We use the definition of $\delta$ to write
$$\delta F_n(\bar{\mu}(\beta,g_1),g_2,\dots,g_{n+1})=$$
$$\sum_{i\neq
1}(-1)^{i+1}\Psi(g_i,F_n(\bar{\mu}(\beta,g_1),g_2,\dots,\hat{g_i},\dots,g_{n+1}))^{\sigma}$$
$$+\Psi(\bar{\mu}(\beta,g_1),F_n(g_2,\dots,g_{n+1}))$$
$$+\sum_{\stackrel{p<q}{p\neq 1}}(-1)^{p+q}F_n(\Phi(g_p,g_q),\bar{\mu}(\beta,g_1),\dots,
\hat{g_p},\dots,\hat{g_q},\dots,g_{n+1}))^{\sigma'}$$
$$+\sum_{q>1}(-1)^{1+q}F_n(\Phi(\bar{\mu}(\beta,g_1),g_q),g_2,\dots,\hat{g_q},\dots,
g_{n+1})^{\sigma''}$$
where $\sigma=(i\,0\,1\,\dots(i-1))$, $\sigma'=(p\,0\,1\dots(p-1))(q\,1\,2\dots(q-1))$ 
and $\sigma''=(q\,2\,3\dots(q-1))$.
We apply the property of $Hom(C,B)$ linearity to $F_n$ and to $\Psi$ in the first three
summands to rewrite the above sum as
$$ \sum_{i\neq 1}(-1)^{i+1}\Psi(g_i,\bar{\nu}(\beta,F_n(g_1,\dots,\hat{g_i},\dots,
g_{n+1}))^{\sigma}$$
$$+\bar{\nu}(\beta,\Psi(g_1,F_n(g_2,\dots,g_{n+1})))$$
$$+\sum_{\stackrel{p<q}{p\neq 1}}(-1)^{p+q}\bar{\nu}(\beta,F_n(\Phi(g_p,g_q),g_1,
\dots,\hat{g_p},\dots,\hat{g_q},\dots,g_{n+1}))^{\sigma'\gamma}$$
$$+\sum_{q>1}(-1)^{1+q}F_n(\Phi(\bar{\mu}(\beta,g_1),g_q),g_2,\dots,\hat{g_q},\dots,
g_{n+1})^{\sigma''}$$
with $\gamma=(2\,0\,1)\in S_{n+2}$.
We next apply the derivation property of $\Psi$ to the first summand and the TD Lie-Rinehart
property to the last summand and obtain
$$\sum_{i\neq 1}(-1)^{i+1}\bar{\nu}(\beta,\Psi(g_i,F_n(g_1,\dots,\hat{g_i},\dots,
g_{n+1}))^{\sigma\tau}$$
$$+\sum_{i\neq 1}(-1)^{i+1}\bar{\nu}(\Psi(g_i,\beta),F_n(g_1,\dots,\hat{g_i},\dots,
g_{n+1}))^{\sigma}$$
$$+\bar{\nu}(\beta,\Psi(g_1,F_n(g_2,\dots,g_{n+1}))$$
$$+\sum_{\stackrel{p<q}{p\neq 1}}(-1)^{p+q}\bar{\nu}(\beta,F_n(\Phi(g_p,g_q),g_1,
\dots,\hat{g_p},\dots,\hat{g_q},\dots,g_{n+1}))^{\sigma'\gamma}$$
$$+\sum_{q>1}(-1)^{1+q}F_n(\bar{\mu}(\beta,\Phi(g_1,g_q)),g_2,\dots,\hat{g_q},\dots,
g_{n+1})^{\sigma''}$$
$$-\sum_{q>1}(-1)^{1+q}F_n(\bar{\mu}(\Psi(g_q,\beta),g_1),g_2,\dots,\hat{g_q},\dots,
g_{n+1})^{\sigma''\gamma'}$$
with $\tau=(0\,1)$ and $\gamma'=(0\,1\,2)$.

Finally, we apply the $Hom(C,B)$ linearity property to the last summand to write it in the form
$$-\sum_{q>1}(-1)^{1+q}\bar{\nu}(\Psi(g_q,\beta),F_n(g_1,\dots,\hat{g_q},\dots,
g_{n+1}))^{\sigma''\gamma'}$$
which will then cancel with the second summand after we note that $\sigma''\gamma'=\sigma$.
In the first summand, we have $\sigma\tau=(1\,2\dots q)$ which we may regard as a cyclic
permutation in $S_{n+1}$.  It then follows that all of the remaining summands yield the
desired $\bar{\nu}(\beta,\delta F_n(g_1,\dots,g_i,\dots,g_{n+1}))$.

For the general case, we first recall that $$\delta F_n\circ(1^{\otimes(i-1)}
\otimes\bar{\mu}\otimes 1^{\otimes(n-i+1)})$$ is induced by the map
$$df_n\circ(1^{\otimes(i-1)}\otimes\mu\otimes 1^{\otimes(n-i+1)})$$ by Proposition \ref{comp}.
We use the skew symmetry of $df_n$ to rewrite this map as $$(-1)^{\sigma}df_n
\circ(\mu\otimes 1^{\otimes n})\circ\sigma\circ\xi$$ where $\xi=(0\,1\,\dots\,(i-1))$ and
$\sigma=(1\,2\,\dots\,i)$.  This map induces the map $$(-1)^{\sigma}(\delta F_n
\circ(\bar{\mu}\otimes 1^{\otimes n})^{\xi\sigma})\circ\sigma\circ\xi$$ which is equal to
$$(-1)^{\sigma}(\bar{\nu}\circ(1\otimes\delta F_n))^{\xi\sigma}
\circ\sigma\circ\xi$$ by the first part of this proof. However, this last map is induced
by the map $$(-1)^{\sigma}(\nu\circ(1\otimes df_n))\circ\sigma\circ\xi$$ which
is equal to $$(-1)^{\sigma}(-1)^{\sigma^{-1}}(\nu\circ(1\otimes df_n))
\circ\sigma^{-1}\circ\sigma\circ\xi$$ by the skew symmetry of $df_n$.  This last map is
equal to $$\nu\circ(1\otimes df_n)\circ\xi$$ which induces the map
$$\bar{\nu}\circ(1\otimes\delta F_n)^{\xi}\circ\xi$$ and the proof is complete.  \qed

\bigskip
\begin{flushleft}
\it E-mail addresses of authors:\\
\sc
gbarnich@ulb.ac.be\\ 
fulp@math.ncsu.edu \\
lada@math.ncsu.edu \\
jds@math.unc.edu\\
\end{flushleft}
\rm

\end{document}